\documentclass[11pt]{amsart}
\usepackage{amssymb}
\usepackage{amsmath}

\newtheorem{thm}{Theorem}[section]
\newtheorem{lem}[thm]{Lemma}
\newtheorem{conj}[thm]{Conjecture}
\newtheorem{prop}[thm]{Proposition}
\newtheorem{cor}[thm]{Corollary}
\theoremstyle{definition}

\newtheorem{exmp}[thm]{Example}
\newtheorem{rem}[thm]{Remark}
\newtheorem{nota}[thm]{Notation}

\newcommand{\beeq}[1]{\begin{eqnarray}\label{#1}}
\newcommand{\eneq}{\end{eqnarray}}

\def\GG{{\mathfrak S}}

\def\spec{{\hbox{\rm{Spec}\,}}}

\def\A{{\mathbb A}}

\def\Q{{\mathbb Q}}
\def\Z{{\mathbb Z}}

\def\P{{\mathbb P}}
\def\C{{\mathbb C}}

\def\O{{\mathcal O}}

\def\H{{\mathfrak H}}
\def\C{{\mathbb C}}

\def\P{{\mathbb P}}

\def\Z{{\mathbb Z}}
\def\Q{{\mathbb Q}}

\def\({\left(}
\def\){\right)}

\def\part{P(n)}

\def\mapl#1{\smash{
\mathop{\longleftarrow}\limits^{#1}}}
\def\mapr#1{\smash{
\mathop{\longrightarrow}\limits^{#1}}}

\def\<{\langle}
\def\>{\rangle}

\begin{document}

\title[On the Motive of the Hilbert scheme]
{On the Motive of the Hilbert scheme of points on a surface}
\keywords{Hilbert scheme, Motives}
\author{Lothar G\"ottsche}
\address{International Center for Theoretical Physics\\Strada Costiera 11\\
34100 Trieste, Italy}
\email{gottsche@ictp.trieste.it}
%\date{\today}

\maketitle\

%\begin{abstract}
%\end{abstract}

%\tableofcontents

\section{Introduction}

Let $S$ be a smooth quasiprojective variety over an algebraically closed field
$k$ of characterstic $0$.
In this note we determine the class $[S^{[n]}]$ of the Hilbert scheme $S^{[n]}$ of subschemes
of length $n$ on $S$ in the Grothendieck ring 
$K_0(V_k)$ of $k$-varieties.
The result expresses $[S^{[n]}]$ in terms of the classes of the symmetric powers
$S^{(l)}=S^l/\GG_l$. Here $\GG_l$ is the symmetric group acting by permutation 
of the factors of $S^l$.

\begin{thm}\label{mainthm}
$$[S^{[n]}]=\sum_{\alpha\in P(n)} [S^{(\alpha)}\times 
\A^{n-|\alpha|}].$$
\end{thm}
Here $P(n)$ is the set of partitions of $n$.
We write a partition $\alpha\in P(n)$ as $(1^{a_1},2^{a_2},\ldots,n^{a_n})$,
where $a_i$ is the number of occurences of $i$ in $\alpha$.
Then the length $|\alpha|$ of $\alpha$ is the sum of the $a_i$
and $S^{(\alpha)}=S^{(a_1)}\times \ldots\times S^{(a_n)}.$

In the case $k=\C$ the cohomology of $S^{[n]}$ has been studied by a number of authors
\cite{E-S1},\cite{Go1},\cite{G-S},\cite{Ch1},\cite{N},\cite{Gr},\cite{L},\cite{dC-M1}.
In particular the Betti numbers and the Hodge structure have been determined.
The class $[X]$ of a smooth projective variety over $\C$ (or 
more generally of a projective variety with finite quotient singularities)
determines its Hodge structure; so Theorem \ref{mainthm} gives a new and 
elementary proof of the corresponding formulas.

By a result of \cite{Gi-Sou} and \cite{Gu-Na}
our result implies that the same formula holds in the Grothendieck ring of 
effective Chow motives. 

Similar arguments apply to the incidence variety 
$S^{[n,n+1]}\subset S^{[n]}\times
S^{[n+1]}$. At the end we give some applications to moduli spaces of 
rank two sheaves on surfaces. 

While I was finishing this paper, the preprint \cite{dC-M2}
appeared, in which the Chow groups and the Chow motive of $S^{[n]}$
 are determined over any field using different methods.

Our approach is mostly motivated by \cite{Ch1} and by \cite{dB}.
Lemma \ref{totaro} plays an important r\^ole in this paper.
I am very thankful to B.~Totaro who proved it for me.
B.~Totaro also pointed out a mistake in an earlier version of the paper
and explained to me   how to  deduce Conjecture \ref{chow} 
from a conjecture of Beilinson and Murre. 
I thank M.S. Narasimhan and K. Paranjape for very 
useful discussions.

\section{Grothendieck rings of varieties and motives and Chow groups}

In this paper let $k$ be an algebraically closed field of characteristic $0$.
Let $K_0(V_k)$ be the Grothendieck ring of $k$-varieties.
This is the abelian group generated by the isomorphism classes of 
$k$-varieties with the relation that $[X\smallsetminus Y]=[X]-[Y]$,
when $Y$ is a closed subvariety of $X$.
The addition and multiplication in this ring are given by
the disjoint union and the product of varieties.

Let $M_k$ be the category of effective Chow motives over $k$.
For the precise definitions and some results about
motives see e.g. \cite{Ma},\cite{Kl},\cite{Sch}.
Let $A^l(X)$ be the $l$-th Chow group of the variety 
$X$ with $\Q$-coefficients. Let $X$, $Y$ be  
smooth projective varieties. Assume $X$ has dimension $d$.
The group $Hom_C(X,Y):=A^d(X\times Y)$ is the group of correspondences  
from $X$ to $Y$ of degree $0$.
An object in $M_k$ is a pair $(X,p)$ where $X$ is a smooth projective 
variety over $k$ and $p\in Hom_C(X,X)$ with $p^2=p$. 
The morphisms are $Hom((X,p),(Y,q))=q Hom_C(X,Y) p$.
There is a contravariant functor $h$ from the category  of smooth projective
varieties to $M_k$ by sending $X$ to $(X,[\Delta_X])$ (where
$\Delta_X\subset X\times X$ is the diagonal) and $f:X\to Y$ to the class of the 
transpose of its graph
$[\Gamma_f']\in A^*(Y\times X)$.
We define $(X,p)\oplus (Y,q):=(X\sqcup Y,p\oplus q)$ 
$(X,p)\otimes (Y,q):=(X\times Y,p\times q)$.

The Grothendieck ring $K_0(M_k)$ is the quotient of the free abelian group 
on the isomorphism classes $[M]$ of effective Chow motives by the subgroup
generated by elements $[M]-[M']-[M'']$ whenever
$M\simeq M'\oplus M''$.
We denote by $[N]$ the class of a motive in $K_0(M_k)$.

\begin{thm}\cite{Gi-Sou},\cite{Gu-Na}.\label{varmot}
Let $k$ be a field of characteristic zero. There exists a unique 
ring homomorphism $\overline h:K_0(V_k)\to K_0(M_k)$ with 
$\overline h([X])=[h(X)]$ for $X$ smooth projective.
\end{thm}

The Lefschetz motive $L$ is defined by $h(\P^1)=1\oplus L$, where $1:=h(pt)$
for $pt=\spec k$.
By $[\P^1]=[pt]+[\A^1]$ we see that $[L]=\overline h([\A^1])$.
So Theorem \ref{mainthm} immediately implies
the identity
$$[h(S^{[n]})]=\sum_{\alpha\in P(n)} \overline h([S^{(\alpha)}])
 [L^{\otimes(n-|\alpha|)}].$$
in  $K_0(M_k)$.

If a finite group $G$ acts on a $k$-variety $X$, the motive
$h(X/G)$ of the quotient can be defined as $(X,\sum_{g\in G} [g])$ where $[g]$ 
is the graph of the action by $g$. Therefore we can associate two a priory different  elements of 
$K_0(M_k)$ to the quotient $X/G$, namely
$[h(X/G)]$ and $\overline h([X/G])$.

\begin{thm}\cite{dB}, Chapter 2, \cite{dB-Na}. $[h(X/G)]=\overline h([X/G])$.
\end{thm}

Therefore we obtain 
\begin{cor}\label{mot}
$$[h(S^{[n]})]=\bigoplus_{\alpha\in P(n)}  [h(S^{(\alpha)})
\otimes  L^{\otimes(n-|\alpha|)}].$$
\end{cor}
 The Chow groups of a Chow motive
$N$ are defined by $A^l(N)=Hom(L^{\otimes l},N)$
and for a smooth projective variety $X$ we have $A^l(h(X))=A^l(X)$.

\begin{rem}\label{chowf}
If $N$ and $M$ are motives with $[N]=[M]$, then there exists a 
motive $P$ with $N\oplus P=M\oplus P$.
By the definition of $\oplus$ and of Chow groups of motives is is evident
that as graded vector spaces 
$A^*(N\oplus P)=A^*(N)\oplus A^*(P)$.
Therefore $A^*(N)\oplus A^*(P)=A^*(M)\oplus A^*(P)$. If the Chow groups
of $P$ are finite dimensional, it follows that $A^*(N)=A^*(M)$.
\end{rem}

We expect that this result holds without the restriction of finite dimensionality.

\begin{conj}\label{chow}
If $N$ and $M$ are effective Chow motives with $[N]=[M]$ in 
$K_0(M_k)$, then 
$M$ and $N$ are isomorphic. In particular they have the same 
Chow groups with rational coefficients. 
\end{conj}
In a previous version of this paper the result of Remark \ref{chowf}
was claimed without the restriction of finite dimensionality. 
The mistake was pointed out to me by B. Totaro.
He also explained to me the following argument how Conjecture \ref{chow}
follows (over any field $k$) from the following conjecture of Beilinson and Murre.

\begin{conj}\label{beil}(see \cite{Ja} Conj. 2.1.).
Let $H^*$ be a Weil cohomology theory. 
For each smooth projective variety $X$ over $k$ and all $j\ge 0$, there exists
a descending filtration $F^\bullet$ on $A^j(X)$ such that
\begin{enumerate}
\item
$F^0 A^j(X)=A^j(X)$ and $F^1A^j(X)$ is the kernel of the cyle map 
$A^j(X)\to H^{2j}(X)$,
\item
$F^rA^i(X)\cdot F^s A^j(X)\subset F^{r+s}A^{i+j}(X)$ for the intersection product,
\item $F^\bullet$ is respected by $f^*$, $f_*$ for morphisms $f:X\to Y$,
\item $F^lA^j(X)=0$ for $l\gg 0$.
\end{enumerate}
\end{conj}
Now we assume Conjecture \ref{beil} and show Conjecture \ref{chow}.
Let $M=(X,p)$ be an effective Chow motive.
Let $R:=End(M)\subset A^*(X\times X)$. The cycle map induces a homomorphism
$R\to End(H^*(X))$. Let $I$ be the kernel.
By the definition of the composition of correspondences and parts 2. and 3. of 
 Conjecture \ref{beil}, we see that for $f\in I$, $f^n\in F^n A^*(X\times X)$.
 So, by part 4., $I$ is nilpotent.
 
Our aim is to show that this implies that effective Chow motives satisfy
 the Krull-Schmidt
Theorem: Every effective Chow motive is the direct sum of finitely many indecomposable
motives, whose isomorphism classes are uniquely defined. This immediately
implies
Conjecture \ref{chow}: If $M, N\in M_k$  with 
$[M]=[N]$, then $M\oplus P\simeq N\oplus P$ for $P\in M_k$.
 By the Krull-Schmidt Theorem it follows that
$M\simeq N$.

In the theorem in Section 3.3 in \cite{Ga-Ro} it is shown that 
an additive category $\mathcal C$ whose isomorphism classes form a set satisfies the
Krull Schmidt Theorem  if the following holds:
Every idempotent in $\mathcal C$ splits and for each object $A$ in 
$\mathcal C$, if we write $R:=End(A)$, then 
$R/rad(R)$ is semisimple and all idempotents in $R/rad(R)$ are the
images of idempotents in $R$. Here $rad(R)$ is the Jacobson radical of $R$.
We check these conditions for the category $M_k$ of effective Chow motives. 
Idempotents split because $M_k$ is pseudoabelian. For a motive $M=(X,p)$
let as above  $R:=End(M)$ and $I=ker(R\to End(H^*(X))$.
Then $I$ is nilpotent and $R/I$ is a finite dimensional $\Q$-algebra.
Since $I$ is nilpotent, $I\subset rad(R)$. So $R/rad(R)$ is a finite
dimensional $\Q$-algebra with radical $0$. So it is semisimple. 
Furthermore we see that $rad(R)$ is nilpotent: $rad(R/I)$
is nilpotent because $R/I$ is finite dimensional. The result follows
for $rad(R)$ because $I$ is nilpotent. 
Then Theorem 1.7.3 in \cite{Be} implies that all idempotents of $R/rad(R)$
lift to idempotents of $R$.

Corollary \ref{mot} and Conjecture  \ref{chow} imply the formulas
\begin{align*}
h(S^{[n]})&=\bigoplus_{\alpha\in P(n)}  h(S^{(\alpha)})
\otimes  L^{\otimes(n-|\alpha|)},\\
A^{i}(S^{[n]})&=\bigoplus_{\alpha\in P(n)} A^{i+|\alpha|-n}(S^{(\alpha)}).
\end{align*}
These formulas have been shown in \cite{dC-M2} over an arbitrary field.

%$
\section{The stratification of $S^{(n)}$ and $S^{[n]}$}
%$
Let $\omega_n:S^{[n]}\to S^{(n)}$ be the Hilbert-Chow morphism, which 
associates
to each
subscheme $Z$ its support with multiplicities.
$S^{[n]}$ and $S^{(n)}$ have a natural stratification parameterized by the
 partitions of $n$,
which has been used before
\cite{Go1},\cite{G-S},\cite{Ch1},\cite{dC-M1} to compute the cohomology
 of $S^{[n]}$. Let $P(n)$ be the set of all partitions of $n$.
 A partition $\alpha=(n_1,\ldots,n_r)$ is also written as 
 $\alpha=(1^{a_1},\ldots,n^{a_n})$, where $a_i$ is the number of occurences of
 $a_i$ in $\alpha$. We write $|\alpha|=r=\sum_i a_i$.
 The corresponding locally closed stratum 
 $S^{(n)}_\alpha$ is the set of all zero-cycles $\xi=n_1x_1+\ldots
 n_r x_r$ with $x_1,\ldots,x_r$ distinct points of $S$.
 We put $S^{[n]}_\alpha:=\omega_n^{-1}(S^{[n]}_\alpha)_{red}.$
 The strata $S^{(n)}_\alpha$ are smooth, but the 
 $S^{[n]}_\alpha$ and the closures  $\overline {S^{(n)}_\alpha}$
usually are singular. 
There is a natural map
$h_\alpha:S^{(\alpha)}\to S^{(n)}; \ (\xi_1,\ldots,\xi_n)\mapsto \sum_{i=1}^n 
i \xi_i$
whose image is the closure $\overline {S^{(n)}_\alpha}$;
in fact it is easy to see that it is the normalization of 
$\overline {S^{(n)}_\alpha}$.
Let $g_\alpha:S^{(\alpha)}\times \A^{n-|\alpha|}\to S^{(n)}$ be the 
composition of the projection to $S^{(\alpha)}$ with 
$h_\alpha$, and let 
$g:\coprod_{\alpha\in P(n)} S^{(\alpha)}\times A^{n-|\alpha|} \to S^{(n)}$
be the map induced by the $g_\alpha$.
A more precise version of Theorem \ref{mainthm} is the following.

\begin{prop} \label{mainprop}
$[g^{-1}(S^{(n)}_\beta)]=[S^{[n]}_\beta]$ in $K_0(V_k)$ for all $\beta\in P(n)$.
\end{prop}

Theorem \ref{mainthm} follows from Proposition \ref{mainprop}
by summing over all $\beta\in P(n)$.

\section{Proof of the main result}

We will determine both sides of the equality in Proposition \ref{mainprop}.
We need some preliminaries.

\begin{rem} 
In the Grothendieck ring of $k$-varieties we have:
\begin{enumerate}
\item If $f:X\to Y$ is a Zariski locally trivial fibre bundle with fibre $F$,
then $[X]=[Y][F]$ (stratify $Y$ such that $f$ is trivial over the strata).
\item If $f:X\to Y$ is a bijective morphism, then 
$[X]=[Y]$. (There is a dense open subset $U\subset X$, on which $f$ is an isomorphism
(here we use $char(k)=0$),
replacing $X$ by $X\smallsetminus U$ and $Y$ by $Y\smallsetminus f(U)$,
we can argue by induction over the dimension.)
\end{enumerate}
\end{rem}

\begin{rem}\label{sym}
If a quasiprojective variety  $X$ has a decomposition 
$X=X_1\sqcup \ldots \sqcup X_l$ into locally
closed subvarieties, then it is immediate to see that 
$\coprod_{n_1+\ldots+n_l=n}\prod_{i=1}^l X_i^{(n_i)}$ is a decomposition of
$X^{(n)}$ into locally closed subvarieties.
Therefore we can define 
$[X]^{(n)}=[X^{(n)}]$ for $X$ a variety, and $([X_1]+\ldots+[X_l])^{(n)}=
\sum_{n_1+\ldots+n_l=n}\prod_{i=1}^l[X_i]^{(n_i)}.$
\end{rem}

\begin{nota}
Let $X$ be a $k$-variety. If $Y$ is a variety with a natural morphism
to $X^{(m)}$ for some $m>0$, we write $Y_*$ for the preimage of the
open subvariety of all zero cycles $\xi\in X^{(m)}$ whose support consists
of $m$ distinct points.
\end{nota}

The following lemma was proved for me by  Burt Totaro.

\begin{lem}\label{totaro}
Let $X$ be a variety over $k$. Let  $p:(X\times \A^l)^{(n)}\to X^{(n)}$
be the obvious  projection.
Then 
$[(X\times \A^l)^{(n)}]=[X^{(n)}\times \A^{nl}]$,
and $[p^{-1}(X^{(n)}_\alpha)]=[X^{(n)}_\alpha\times \A^{nl}]$ for all 
$\alpha\in P(n)$.
\end{lem}

\begin{proof}
The first statement follows from the second.
It is enough to treat the case $l=1$;  the 
general case follows by trivial induction. 
There is a cartesian diagram
$$\begin{matrix}
 X^{|\alpha|}_*\times \prod_{i=1}^n((\A^1)^{(i)})^{a_i}&\mapr{\overline p} &
 X^{|\alpha|}_*\\
\downarrow&&\downarrow\\
p^{-1}(X^{(n)}_\alpha)&\mapr{p}&X^{(n)}_\alpha
\end{matrix}
$$
By the fundamental theorem of symmetric functions the fibre product is
$X^{|\alpha|}\times \prod_i (\A^i)^{a_i}=X^{|\alpha|}_*\times \A^n.$
So we get an \'etale trivialization of $p$. Any two trivializations are 
related by the action of the group 
$\GG_{a_1}\times \ldots\times \GG_{a_n}$ by
reordering the factors in the $X^{a_i}$.
This acts on the fibres of $\overline p$  by reordering the factors $(\A^1)^{(i)}$
in the $((\A^1){(i)})^{a_i}$.
Choosing an origin in $\A^1$ determines an origin in $\A^n$, and 
the action becomes linear on on $k^n$, so $p^{-1}(X^{(n)}_\alpha)$
is an \'etale locally trivial vector bundle over $X^{(n)}_\alpha$.
Therefore, by Hilbert Theorem 90 \cite{Se} p. 1.24, it is locally trivial 
in the Zariski topology.
Thus we get 
$[p^{-1}(X^{(n)}_\alpha]=[X^{(n)}_\alpha\times \A^n]$.
\end{proof}

Now we determine the right hand side in Proposition \ref{mainprop}.
Let $R_n=Hilb^n(\A^2,0)$ be the punctual Hilbert scheme of subschemes of
length $n$ of $\A^2$ concentrated in $0$.
Then by \cite{E-S1} $R_n$ has a cell decomposition and
$$[R_n]=\sum_{\beta\in P(n)}[\A^{n-|\beta|}].$$

\begin{lem}\label{hilquot}
$[S^{[n]}_\alpha]= \big[\bigl(\prod_{i=1}^n(S\times R_i)^{(a_i)}\bigr)_*].$
\end{lem}

\begin{proof}
By Lemma 2.1.4 of \cite{Go2} $S^{[l]}_{(l)}$ is a locally trivial fibre 
bundle over $S$ with  fibre $R_l$,
thus $[S^{[l]}_{(l)}]=[S\times R_l]$.
There is a natural morphism 
$f:\bigl(\prod_{i=1}^n (S^{[i]}_{(i)})^{(a_i)}\bigr)_*
\to S^{[n]}_\alpha$ defined on $T$-valued points by sending 
$(Z_1,\ldots,Z_n)$ to $\coprod_{i=1}^n
Z_i$.
$f$ is obviously invariant under the action  of $\GG_{a_1}\times \ldots\times
\GG_{a_n}$  by permuting the factors in the $(S^{[i]}_{(i)})^{(a_i)}$,
and  the induced morphism from the quotient to  $S^{[n]}_\alpha$ 
induces a bijection on $k$-valued points. 
This implies 
$[S^{[n]}_\alpha]=\big[\bigl(\prod_{i=1}^n 
(S^{[i]}_{(i)})^{(a_i)}\bigr)_*\big]$.
\end{proof}

\begin{nota}
\begin{enumerate}
\item
For any $x\in S$ and any $\xi\in S^{(n)}$ we call $m_x(\xi)$ the multiplicity with 
which $x$ occurs in $\xi$.
\item
We denote $P:=\bigcup_{n>0} P(n)$. For $\alpha=(1^{a_1},2^{a_2},\ldots)\in P(n)$, 
$\beta:=(1^{b_1},2^{b_2},\ldots)\in P(m)$ and $l\in \Z_{\ge 0}$ we denote
$l\alpha:=(1^{la_1},2^{la_2}\ldots)\in P(nl)$, 
$\alpha+\beta:=(1^{a_1+b_1},2^{a_2+b_2}\ldots)\in P(n+m)$.
\end{enumerate}
\end{nota}

\begin{lem}\label{func}
$[S^{[n]}_\alpha]=\big[
 \coprod_{f}\bigl(\prod_{\beta\in P} S^{(f(\beta))}\bigr)_*\times
\A^{n-\sum_{\beta\in P} f(\beta)|\beta|}\big].$

Here $f$ runs through the  $f:P\to \Z_{\ge 0}$
with $\sum_{\beta\in P(i)} f(\beta)=a_i$ for all $i$.
\end{lem}

\begin{proof}
By Lemma \ref{hilquot} we get
$[S^{[n]}_\alpha]= \big[\bigr(\prod_{i=1}^n\big( \coprod_{\beta_i\in P(i)}
S\times \A^{i-|\beta_i|}\big)^{(a_i)}\bigl)_*\big].$
By Remark \ref{sym} this implies
$[S^{[n]}_\alpha]=\big[\bigr(\prod_{i=1}^n\coprod_{f_i}\prod_{\beta_i\in P(i)}
 (S\times \A^{i-|\beta_i|})^{(f_i(\beta_i)}\bigl)_*\big],$
where the $f_i$ rund through the $f_i:P(i)\to \Z_{\ge0}$ with
$\sum_{\beta\in P(i)} f_i(\beta)=a_i$.
The result follows by Lemma \ref{totaro}.
\end{proof}

Now we determine the left hand side in Proposition \ref{mainprop}.
Let $\beta=(1^{b_1},2^{b_2},\ldots)=(n_1,\ldots,n_r)\in P(n)$.

\begin{lem}\label{ofunc}
$h_\alpha^{-1}S^{(n)}_\beta\simeq  \coprod_{f} 
\bigl(\prod_{\gamma\in P} 
S^{(f(\gamma))}\bigr)_*.$

Here the sum is over all functions $f:P\to \Z_{\ge 0}$ with 
$\sum_{\gamma\in P(i)} f(\gamma)=b_i$ and
$\sum_{\gamma\in P} f(\gamma)\gamma =\alpha$.
\end{lem}

Using that $g_\alpha^{-1}(S^{(n)}_\beta)=
h_\alpha^{-1}(S^{(n)}_\beta)\times \A^{n-|\alpha|}$,
Proposition \ref{mainprop} follows immediately from Lemma \ref{func} and
Lemma \ref{ofunc} by summing over all $\beta$.

\begin{proof}
Any $\xi=(\xi_1,\ldots,\xi_n)\in h_\alpha^{-1}S^{(n)}_\beta$ induces a map 
$f_\xi:P\to \Z_{\ge 0}$ as follows.
Let $h_\alpha(\xi)=\sum_{i=1}^r n_i x_i$. For all $x\in S$ let 
$\gamma_x(\xi):=\sum_{j=1}^n (j^{m_{x}(\xi_j)})$. For each $i=1,\ldots r$ 
we have
$\sum_{j=1}^n j m_{x_i}(\xi_j)=n_i$, so $\gamma_{x_i}(\xi)\in P(n_i)$;
furthermore $\sum_i \gamma_{x_i}(\xi)=\alpha$. We define $f_\xi:P\to\Z_{\ge 0}$ by 
$f_\xi(\gamma):=\# \big\{x \in S\bigm | \gamma_x(\xi)=\gamma\}.$ 
Then $\sum_{\gamma\in P(j)} f_\xi(\gamma)=b_j$ and
$\sum_{\gamma\in P}f_\xi(\gamma)\gamma=\gamma_{x_1}(\xi) +\ldots
+\gamma_{x_r}(\xi)=\alpha.$

Now fix $f:P\to \Z_{\ge 0}$ with the above properties. 
Let $S^{(\alpha)}_{f}:=\big\{\xi\in S^{(\alpha)}\bigm| 
f_\xi=f\}$. We claim that 
$S^{(\alpha)}_{f}\simeq 
\bigl(\prod_{\gamma\in P}S^{(f(\gamma))}\bigr)_*.$

For $\xi\in S^{(\alpha)}_\Gamma$ define 
$\phi(\xi)=(\phi(\xi)_\gamma)_{\gamma\in P}$ by letting
$\phi(\xi)_\gamma\in S^{(f(\gamma))}$ be the sum over all $x\in S$ with 
$\gamma_x(\xi)=\gamma$.
For $\zeta=(\zeta_\gamma)_{\gamma\in P}$ with $\zeta_\gamma\in S^{(f(\gamma))}$
let $\psi(\zeta):=
(\xi_1,\ldots,\xi_n)$ with 
$\xi_i=\sum_{\gamma\in P} c_i\zeta_{\gamma}\in 
S^{(a_i)},$ where we write $\gamma=(1^{c_1},2^{c_2},\ldots)$. 
It is straightforward from the definitions  that $\phi$ and $\psi$ are 
inverse to each other.
\end{proof}

\begin{exmp} \label{rathilb}
\begin{enumerate}
\item Let $S$ be a projective rational surface. 
Then $[S]=[\A^0]+b[\A^1]+[\A^2]$ for suitable $b>0$, and
$$
\sum_{n\ge 0}[S^{[n]}]t^n=\prod_{l>0} \frac{1}{(1-[\A^{l-1}]t^l)
(1-[\A^{l}]t^l)^b(1-[\A^{l+1}]t^l)}.
$$
\item Let $S$ be a ruled surface over a curve $C$.
Then 
$$\sum_{n\ge 0}[S^{[n]}]t^n=
\prod_{l>0}\Big(\sum_{m\ge 0} [C^{(m)}\times  \A^{m(l-1)}]t^{ml}\Big)
\Big(\sum_{m\ge 0} [C^{(m)}\times  \A^{ml}]t^{ml}\Big).
$$
\item If $\widehat S$ is the blowup of $S$ in a point then
$$\sum_{n\ge 0}[\widehat S^{[n]}]t^n=\frac{\sum_{n\ge 0}[S^{[n]}]t^n}
{\prod_{l>0}(1-[\A^l]t^l)}.$$
\end{enumerate}
This follows from Theorem \ref{mainthm} and Lemma \ref{totaro}.
\end{exmp}

\section{The incidence variety}
Similar but simpler arguments to those for $S^{[n]}$ can be used for
 the incidence variety
$S^{[n,n+1]}:=\big\{(Z,W)\in S^{[n]}\times S^{[n+1]}\bigm|
Z\subset W\}$, which plays a r\^ole in inductive arguments for $S^{[n]}$
\cite{E-S2},\cite{E-G-L}. The Hodge numbers of $S^{[n,n+1]}$ were computed
in \cite{Ch1}.

\begin{thm}\label{incid}
$\displaystyle{[S^{[n,n+1]}]=\sum_{l=0}^n [S\times S^{[l]}\times \A^{n-l}]}.$
\end{thm}

By Theorem \ref{varmot} this immediately implies

\begin{cor}
$\displaystyle{[ h(S^{[n,n+1]})]=
\bigoplus_{l=0}^n [h(S\times S^{[l]})\otimes L^{\otimes(n-l)}]}.$
\end{cor}

For the proof we introduce a stratification of $S^{[n,n+1]}$.
Let $$\overline \omega:S^{[n,n+1]}\to S\times S^{(n)},(Z,W)\mapsto 
(\omega_{n+1}(W)-\omega_n(Z),\omega_n(Z)).$$
For $0\le m\le n$ let 
$(S\times S^{(n)})_m:=\big\{(x,\xi)\in S\times S^{(n)}\bigm| m_x(\xi)=m\big\}$,
and let 
$S^{[n,n+1]}_m:=\overline \omega^{-1}((S\times S^{(n)})_m)_{red}$.
The  $(S\times S^{(n)})_m$ and the $S^{[n,n+1]}_m$ form  stratifications
of $S\times S^{(n)}$ and  $S^{[n,n+1]}$ respectively
into locally closed subvarieties. 
Let 
$$\overline g:\coprod_{m=0}^n S\times S^{[m]}\times \A^{n-m}
\to S\times S^{(n)}; (x,Z,a)\mapsto (x, (n-m)x +\omega_m(Z)).$$
Then Theorem \ref{incid} follows from the following.

\begin{prop}\label{incidprop}
$\displaystyle{[S^{[n,n+1]}_m]=[\overline g^{-1}((S\times S^{(n)})_m).]}$
\end{prop}

\begin{proof}
If $X$ is a variety with a natural map to $S\times S^{(m)}$ for
some $m\ge 0$, we will 
write $X_0$ for the preimage of the locus of $(x,\xi)$ with 
$x\not \in supp(\xi)$.
Let $(\A^2,0)^{[n,n+1]}:=\big\{(Z,W)\in (\A^2)^{[n,n+1]}\bigm| Z\subset W,
supp(W)=\{0\}\big\}$ with the reduced structure. In \cite{Ch2} it is shown
that $(\A^2,0)^{[n,n+1]}$ has a cell decomposition. Her formula for 
the numbers of cells of different dimensions implies that 
$[(\A^2,0)^{[n,n+1]}]=\sum_{l=0}^n [R_l\times \A^{n-l}].$
We now determine $[S^{[n,n+1]}_m]$.
First it is easy to see analoguously to the case of $S^{[n]}_{(n)}$ 
in \cite{Go2} that $S^{[n,n+1]}_n$ is a locally trivial fibre bundle over 
$S$ with fibre $(\A^2,0)^{[n,n+1]}$.
Therefore 
$$[S^{[n,n+1]}_n]=\sum_{l=0}^n [S\times R_l\times \A^{n-l}]=
\sum_{l=0}^n [S^{[l]}_{(l)}\times \A^{n-l}].$$
There is a natural morphism 
$\sigma:(S^{[m,m+1]}_m\times S^{[n-m]})_0\to S^{[n,n+1]}_m$
given on $T$-valued points by $((Z,W),X)\mapsto (Z\sqcup X,W\sqcup X)$.
$\sigma$ is obviously a bijection on $k$ valued points. Thus we get
$$[S^{[n,n+1]}_m]=[(S^{[m,m+1]}_m\times S^{[m-n]})_0]=
\sum_{l=0}^m[(S^{[l]}_{(l)}\times S^{[n-m]})_0\times \A^{m-l}].$$
Now we determine $\overline g^{-1}(S\times S^{(n)})_m)$.
Let $(S\times S^{[m]})_l=\big\{(x,Z)\in S\times S^{[m]}\bigm| 
len_x(Z)=l\big\}.$
Then $\overline g^{-1}((S\times S^{(n)})_m)=\coprod_{l=0}^m 
(S\times S^{[n-m+l]})_l\times \A^{m-l}$.
Furthermore we have a morphism
$\phi:(S^{[l]}_{(l)}\times S^{[m-l]})_0\to (S\times S^{[m]})_l$
sending $(Y,Z)$ to $(supp(Y),Y\sqcup Z)$, 
which is bijective on $k$-valued points.
Thus $[g^{-1}((S\times S^{(n)})_m)]=
\sum_{l=0}^m[(S^{[l]}_{(l)}\times S^{[n-m]})_0\times \A^{m-l}].$

\end{proof}

\section{Moduli of stable sheaves}

Let $S$ be a projective surface over $k$. Fix $C\in Pic(S)$ and
let $H$ be an ample line bundle on $S$. We denote by $NS(S)$ the
Picard group of $S$ modulo numerical equivalence.
The moduli space $M_S^H(C,d)$ of $H$-semistable rank $2$ torsion-free
sheaves $E$ with $det(E)=C$ and $c_2-C^2/4=d$ depends on $H$ via a system 
of walls and
chambers. This dependence has been studied and used by various authors
(e.g. \cite{Q1},\cite{F-Q},\cite{E-G},\cite{Go3}). 
A class $\xi\in NS(S)+C/2$ is of 
type $(C,d)$ if $0<-\xi^2\le d$ and there exists an ample divisor $H$
with $\xi\cdot H=0$. In this case we say that $H$ lies on the  
corresponding wall. Ample divisors $H, L\in Pic(S)$ are separated by $\xi$ if
$(\xi\cdot H)(\xi\cdot L)<0$.
Assume that $H$ and $L$ do not lie on a wall of type $(C,d)$.
If they are not separated by a class of type $(C,d)$ we say that they lie in 
the same chamber of type $(C,d)$. In this case $M_S^L(C,d)=M_S^H(C,d)$.
More generally let $F\in Pic(S)$ be nef with 
$F^2\ge 0$ and  $FC$ odd.
Then $F+nH$ is ample for any ample divisor $H$,
and
for $n$ sufficiently large the chamber of $F+nH$ does not depend on $H$.
We will write $M_S^F(C,d):=M_S^{F+nH}(C,d)$ for $n\gg 0$.
Let $L$, $H$ be ample divisors not on a wall of type $(C,d)$.
If $2\xi+K_S$ is not 
effective for all $\xi$ separating $L$ and $H$ (we say that $L$ and $H$ are separated
only by good walls), then 
$M_S^L(C,d)$ is obtained from $M_S^H(C,d)$ by successively blowing up 
along projective space bundles over products $S^{[l]}\times S^{[n]}$
of Hilbert schemes of points followed by blowdowns of the exceptional divisor 
to another projective space bundle over $S^{[l]}\times S^{[n]}$.
Therefore the proof of Theorem 3.4 in \cite{Go3} shows:

\begin{prop}\label{wallcr}
Let $L$, $H$ be ample divisors not on a wall of type $(C,d)$ 
separated only by good walls.
Then
\begin{align*}[M_S^L&(C,d)]-[M_S^H(C,d)]= [Pic^0(S)]\cdot \\
&\sum_{\xi} [(S\sqcup S)^{[d+\xi^2]}]
\big([\P^{d-\xi^2+\xi K_S-\chi(\O_S)-1}]-
[\P^{d-\xi^2-\xi K_S-\chi(\O_S)-1}]\big).
\end{align*}
The sum is over all $\xi$ of type $(C,d)$  with 
$\xi L>0>\xi H$, and we use the convention $\P^{-1}=\emptyset$.
\end{prop}

\begin{cor}
Under the assumptions of Proposition \ref{wallcr}, if $S$ is  
a rational surface, then 
$[M_S^L(C,d)]-[M_S^H(C,d)]$ is a $\Z$-linear combination of 
the $[\A^l]$ with $l\le 4d-3$.
\end{cor}
This follows from Proposition \ref{wallcr}, and Example
\ref{rathilb}.

\begin{cor}
If $K_S$ is numerically trivial, then
$[M_S^H(C,d)]$ does not depend on $H$ as long as $H$ does not lie on 
a wall.
\end{cor}

In \cite{G-Z} and \cite{Go4} Theta functions for indefinite lattices
were introduced to study the wallcrossing. 
Let $\Gamma$ be $NS(S)$ with the negative of the intersection form
as quadratic form, which we denote by $\<\ ,\ \>$.
Then for $F,G\in Pic(S)$ with $F^2\ge 0$, $G^2\ge 0$, $F\cdot G>0$, we define
$$\Theta^{F,G}_{\Gamma,C}(\tau,x)
:=\sum_{\xi\in \Gamma+C/2}( \mu(\<\xi,F\>-\mu(\<\xi,G\>)q^{\<\xi,\xi\>/2} 
e^{2\pi i\<\xi,x\>}.$$
Here $\mu(t)=1$ of $t\ge 0$, and $\mu(t)=0$ otherwise, and 
$q=e^{2\pi i\tau}$ for $\tau$ in the complex upper half plane $\H$ and
$x\in \Gamma\otimes_\Z\C$. This function is
defined on a suitable open subset of $\H\times \Gamma\otimes_\Z\C$
and has a meromorphic extension to the whole of 
$\Gamma\otimes_\Z\C$. 
For a meromorphic function $f:\H\times \Gamma\otimes_\Z \C$ and $v\in
\Gamma\otimes_Z \Q$ we write
$$f|_v(\tau,x):=q^{\<v,v\>/2}e^{2\pi i \<v,x\>}f(\tau,x+v\tau).$$
Then $\Theta^{F,G}_{\Gamma,C}(\tau,x)=\Theta^{F,G}_{\Gamma}|_{C/2}(\tau,x)$,
where we have written $\Theta^{F,G}_{\Gamma}:=\Theta^{F,G}_{\Gamma,0}.$

The reason for introducing these theta functions was that they
can be expressed in terms of standard theta functions in case
$F^2=G^2=0$. In the rest of this section we
write $y:=e^{2\pi i z}$ for $z$ a complex variable. 
Recall the standard theta functions
$$\Theta_{\mu,\nu}(\tau,z):=\sum_{n\in Z}(-1)^{n\nu}
q^{(n+\mu/2)^2/2}y^{n+\mu/2}
\quad (\mu,\nu\in \{0,1\}).$$ 
If
$F^2>0$ and $G^2>0$ or $F\cdot C$ and $G\cdot C$ are odd, then
for every $L\in \Gamma$, $\Theta^{F,G}_{\Gamma,C}(\tau,Lz)$
is a power series in $q^{1/8}$ with coefficients Laurent polynomials 
in $y^{1/2}$.
We will write 
$\Theta^{F,G}_{\Gamma,C}(2\tau,K_Sz)^*$ for the power series in $t^{1/4}$
with coefficients Laurent polynomials in $[\A^1]$, which we obtain by 
replacing $y$ by $[\A^1]$ and $q$ by $[\A^2]t$ in  
$\Theta^{F,G}_{\Gamma,C}(2\tau,K_Sz)$. There are no half integer powers
of $[\A^1]$, because $K_S$ is characteristic.

The following follows from Theorem 3.4 in \cite{Go3} in the same way as
Theorem 4.1 in \cite{Go4}.

\begin{cor}\label{wallthet}
Assume $H,L$ do not lie on a wall of type $(C,d)$ for any $d$,
and are separated only by good walls. Then, in $K_0(V_k)$
\begin{align*}
\sum_{d\ge 0}&([M_S^H(C,d)]-[M_S^L(C,d)])t^d&\\
&=[Pic^0(S)]\Big(\sum_{n\ge 0} [S^{[n]}][\A^n]t^n\Big)^2
\frac{\Theta^{L,H}_{\Gamma,C}(2\tau,K_S z)^*}{
[\A^1]^{-\chi(\O_S)}([\A^1]-1)}.
\end{align*}
On the r.h.s. we mean that, after multiplying out, all
negative powers of $[\A^1]$ vanish.
\end{cor}

In \cite{L-Q1}, \cite{L-Q2} a blowup formula was proven for the 
Euler numbers and the virtual
Hodge numbers of moduli spaces of stable rank $2$ sheaves on surfaces.
In \cite{Ka} a blowup formula is proven for principal bundles. 
We can show that the formula of  \cite{L-Q1}, \cite{L-Q2} holds in
 $K_0(V_k)$.
Let $H$ be ample on $S$ and assume that $C\cdot H$ is odd. 
Let $\widehat S$ be the blowup of $S$ in a point and denote by $E$
the exceptional divisor, we denote by $H$ also the pullback of $H$ to 
$\widehat S$. 

\begin{thm}\label{blowup} Assume $k=\C$.
Let $a\in \{0,1\}$ then
\begin{align*}
\frac{\sum_{d\ge 0} [M_{\widehat S}^H(C+aE,d)]t^{d}}
{\sum_{d\ge 0} [M_{ S}^H(C,d)]t^{d}}
&=
\frac{\sum_{n\in \Z} [\A^{\binom{2n+a+1}{2}}]t^{(n+\frac{a}{2})^2}}
{\prod_{l>0}(1-[\A^{2l}]t^{l})}.
\end{align*}
\end{thm}

\begin{proof}
In \cite{L-Q1} the authors use virtual Hodge polynomials
$e(X:x,y)$ in order
to show that there exists a universal power series
$Z_a(x,y,t)$ such that 
$$\sum_{d\ge 0} e(M_{\widehat S}^H(C+aE,d),:x,y)t^{d}=Z_a(x,y,z)
\Big(\sum_{d\ge 0} e(M_{ S}^H(C,d):x,y)t^{d}\Big).$$
To do this they only use the  basic property of virtual Hodge polynomials
that $e(X\smallsetminus Y:x,y)=e(X:x,y)-e(Y:x,y)$ for $Y$ a closed subvariety 
of $X$.
So their proof shows that there is a universal power series $Y_a(t)$ in 
$K_0(V_k)[[t]]$, such that 
$$\sum_{d\ge 0} [M_{\widehat S}^H(C+aE,d)]t^{d}=
Y_a(t)\Big(\sum_{d\ge 0} [M_{ S}^H(C,d)]t^{d}\Big).$$
In the paper \cite{L-Q2} they compare the wallcrossing on a rational
ruled surface and its blowup in a point to determine $Z_a(x,y,t)$.
We can translate their argument into our language, where it proves the
theorem:
Let $H$, $L$ be ample on $S$ with $C\cdot H$ and $C\cdot L$ odd.
Write $M_{d}:=[M_{S}^H(C,d)]-[M_{S}^L(C,d)]$ and 
$M_{a,d}:=[M_{\widehat S}^H(C+aE,d)]-[M_{\widehat S}^L(C+aE,d)]$.
Write $\Gamma=H^2(S,\Z)$.
Then, as $Y_a(t)$ is universal, 
we get, using Corollary \ref{wallthet},
\begin{align*}
Y_a(t)&=\frac{\sum_{d\ge 0} M_{a,d} t^d}{\sum_{d\ge 0} M_{d}t^d}
=\frac{
\Theta^{L,H}_{\Gamma\oplus\<E\>,C+aE}(2\tau,K_{\widehat S}z)^*
\Big(\sum_{n\ge 0}[\widehat S^{[n]}]t^n\Big)^2 }
{\Theta^{L,H}_{\Gamma,C}(2\tau,K_Sz)^*
\Big(\sum_{n\ge 0} [S^{[n]}]t^n\Big)^2 }
\end{align*}
By definition  it is obvious that
$\Theta^{L,H}_{\Gamma\oplus\<E\>,C+aE}(\tau,K_{\widehat S}z)=
\theta_{a,0}(\tau,z)
\Theta^{L,H}_{\Gamma,C}(\tau,K_Sz)$. 
The result follows by Example \ref{rathilb}.(3) and the identity
$$\theta_{a,0}(2\tau,z)^*=
\sum_{n\in \Z} [\A^{\binom{2n+a+1}{2}}]t^{(n+\frac{a}{2})^2}.$$
\end{proof}

Corollary \ref{wallthet} and Theorem \ref{blowup}
imply in particular that the computations of \cite{Go4} hold in
the Grothendieck group of varieties (replacing $y$ in the formulas
there by $[\A^1]$).
As there, we use the following elementary fact.

\begin{rem} \label{vanish}
Let $S$ be the blowup of a ruled surface in finitely 
many points
and $F$ the pullback of a fibre of the ruling. Assume $F\cdot C$ is odd. 
Then $M_S^F(C,d)=\emptyset$ for all $d$.
\end{rem}

 In particular we get the following results.

\begin{cor} \label{ratcor} Let $S$ be a rational surface.
\begin{enumerate}
\item Let $C\in Pic(S)\setminus \{0\}$.
Let $H$ be an ample divisor not on a wall of type $(C,d)$
and assume that $K_S\cdot H\le 0$, then
$[M_S^H(C,d)]$ is a $\Z$-linear combination of $[\A^l]$ with $l\le 4d-3$.
(\cite{Go4}, Prop.4.9),
\item Let $S$ be the blowup of $\P^2$ in 9 points, with exceptional 
divisors $E_1,\ldots,E_9$ and let $H$ be the pullback of the 
hyperplane class. Let $F:=3H-E_1-\ldots-E_9$, and let 
$C\in H^2(S,\Z)$ with $C^2$ odd. 

Then
$[M_S^F(C,d)]=[S^{[2d-3/2]}]$.
(\cite{Go4}, Thm.7.3).
\end{enumerate}
\end{cor}

\begin{rem}
For $S$ a rational surface it follows from \cite{dC-M2} that
the rational Chow groups of the $S^{[n]}$ are finite dimensional.
By Remark \ref{chowf} and the discussion before Proposition \ref{wallcr}
it follows that under the conditions of Corollary \ref{ratcor}.1.
the calculations in \cite{Go4} hold also in the Chow ring of
$M_S^H(C,d)$. In particular in 2. we get that $M_S^F(C,d)$ and
$S^{[2d-3/2]}$ have the same Chow groups with rational coefficients. 
\end{rem}

As a final example we determine the classes of some moduli
spaces over a ruled surface $S$ over an elliptic curve $E$
with a section $\sigma$ with minimal self-intersection $\sigma^2=1$.

\begin{prop}
Let $F$ be the class of a fibre of the ruling.
We write $G:=2\sigma-F$.
Let $C\in NS(S)$ with $C^2$ odd.
Then $[M_S^G(C,d)]$ is the coefficient of $t^{2d-1/2}$ in
\begin{align*}
[E]&\prod_{m>0}\prod_{i=0,1}\Bigg(1+[E]
\Bigg(\sum_{l\ge 1}[\P^{l-1}][\A^{m(2l-i)}]t^{2lm}\Bigg)\Bigg)^2\cdot\\
&\prod_{n>0} ((1-[\A^{n-1}]t^n)(1-[\A^{n}]t^n)^2(1-[\A^{n+1}]t^n))^{(-1)^n}
\end{align*}
\end{prop}

\begin{proof}
$C^2$ odd implies $C\cdot F$ odd and $C\cdot G$ odd,
therefore Remark \ref{vanish} implies that $M_S^F(C,d)=\emptyset$.
By Proposition \ref{wallcr} $[M_S^G(C,d)]$ depends only on the numerical 
equivalence class of $C$. Therefore we can assume that
$C=\sigma=\frac{G+F}{2}$ of $C=\sigma-F=\frac{G-F}{2}$.
We write $\Gamma=NS(S)$ with the negative of the intersection form.
We note that $K_S=-G$.
By Corollary \ref{wallthet} $[M_S^G(C,d)]$ is the coefficient of $t^d$ in
$$[E]
\Big(\sum_{n\ge 0}[S^{[n]}][\A^n]t^n\Big)^2
\frac{\Theta^{F,G}_{\Gamma,\frac{F+G}{2}}(2\tau,-Gx)^*+
\Theta^{F,G}_{\Gamma,\frac{F-G}{2}}(2\tau,-Gx)^*}{[\A^1]-1}.
$$
F or $n\ge 1$ $E^{(n)}$ is  a locally trivial bundle over $E$ with fibre  
$\P^{n-1}$, thus $[E^{(n)}]=[E\times \P^{n-1}]$.
Using Example \ref{rathilb},
we only need to show that
\begin{align*}
\Theta^{F,G}_{\Gamma,\frac{F+G}{2}}&(2\tau,-Gz)+
\Theta^{F,G}_{\Gamma,\frac{F-G}{2}}(2\tau,-Gz)\\
&=(y^{1/2}-y^{-1/2})
q^{1/4}\prod_{n>0}\big((1-q^{n/2}y^{-1})(1-q^{n/2})^2(1-q^{n/2}y)\big)^{(-1)^n}.
\end{align*}
Let $L$ be the lattice generated by $F/2,G/2$. 
$0$ and $G/2$ are a basis of $L$ modulo $\Gamma$.
Therefore $\Theta^{F,G}_{\Gamma}+\Theta^{F,G}_{\Gamma,G}=\Theta^{F,G}_{L}$
and $\Theta^{F,G}_{\Gamma,\frac{F+G}{2}}+
\Theta^{F,G}_{\Gamma,\frac{F-G}{2}}=\Theta^{F,G}_{L,\frac{F+G}{2}}$.
By $F^2=G^2=0$ and  $F\cdot G=2$, 
formula (2.14) in \cite{Go4} gives that
$$\Theta_{L}^{F,G}(2\tau,x)=
\frac{\eta(\tau)^3\theta_{1,1}(\tau,\<(-F+G)/2,x\>)}
{\theta_{1,1}(\tau,\<-F/2,x\>)\theta_{1,1}(\tau,\<G/2,x\>)}.$$
Easy computations give 
$$
\Theta_{L,\frac{F+G}{2}}^{F,G}(2\tau,x)=
\Theta_{L}^{F,G}|_{\frac{F+G}{2}}(2\tau,x)=
\frac{\eta(\tau)^3\theta_{1,1}(\tau,\<(-F+G)/2,x\>)}
{\theta_{0,1}(\tau,\<-F/2,x\>)\theta_{0,1}(\tau,\<G/2,x\>)},$$
$$\Theta_{L,\frac{F+G}{2}}^{F,G}(2\tau,-Gz)=
\frac{\eta(\tau)^3\theta_{1,1}(\tau,z)}
{\theta_{0,1}(\tau,z)\theta_{0,1}(\tau,0)}.
$$
The result now follows from the product formulas
\begin{align*}
\theta_{1,1}(\tau,z)&=q^{\frac{1}{8}}(y^{\frac{1}{2}}-y^{-\frac{1}{2}})
\prod_{n>0}(1-q^n)(1-q^{n}y)((1-q^{n}y^{-1}),\\
\nonumber\theta_{0,1}(\tau,z)&=
\prod_{n>0}(1-q^n)(1-q^{n-\frac{1}{2}}y)((1-q^{n-\frac{1}{2}}y^{-1}).
\end{align*}
\end{proof}

\section{Remarks and Speculations}

Recently there has been a lot of interest in motivic integration
\cite{Ko}, \cite{D-L}, \cite{Lo}. This is a method to determine 
the class of a variety in  a ring $\widehat M_k$ 
which is obtained from $K_0(V_k)$ via localization
at the class of $[\A^1]$ and a suitable completion.
In many cases the result of motivic integration is an 
identity between the classes of smooth projective varieties or 
projective varieties with finite quotient singularities 
in $\widehat M_k$,
and one would expect that the identity does indeed hold in $K_0(V_k)$.
Using Conjecture \ref{chow} this conjecturally gives an
 isomorphism of motives and of the Chow groups 
with rational coefficients.

Two important instances of this are the following:

(1)  Let $X,Y$ be a smooth projective birational $k$-varieties with 
$K_X=K_Y=0$.
Then   motivic integration is used to show that 
$X$ and $Y$ have the same Hodge numbers (\cite{Ko},\cite{D-L}) (that they have the same
Betti numbers was shown before in \cite{Ba} via $p$-adic integration). 
In fact  they have
the same class in $\widehat M_k$. It should be true that $[X]=[Y]$ in 
$K_0(V_k)$.

This happens e.g.  in some cases for moduli spaces of K3-surfaces:
Let $S$ be a K3 surface with $Pic(S)=\Z L$ for $L$ an ample divisor.
Then $[M_S^L(L,L^2/2+3)]=[S^{[L^2/2+3]}]$.
This follows from the proof of Proposition 1.9 in \cite{G-H}
We write $M:=M_S^L(L,L^2/2+3)$, $X:=S^{[L^2/2+3]}$.
There is a  diagram of birational maps
$M\mapl{\phi} N\mapr{\psi} X$ together with stratifications
$M=\coprod M_l$, $N=\coprod N_l$, $X=\coprod X_l$ with
$N_l=\phi^{-1}M_k=\psi^{-1}X_l$, such that 
$N_l\to M_l$ and $N_l\to X_l$ are $\P^{l-1}$-bundles.

(2) In a similar way the results of motivic integration on the
McKay correspondence \cite{D-L},\cite{R}
make it seem likely that the following holds in the Grothendieck ring 
of varieties.
Let $X$ be a smooth projective variety  acted upon by a finite group 
$G\subset Sl(n,\C)$ such that the action preserves the canonical divisor.
Let $Y$ be a crepant resolution of $X/G$.
Choosing an eigenbasis, we can write $g=diag(\epsilon^{a_1},\ldots,
\epsilon^{a_n})$, where $\epsilon$ is a primitive $r^{th}$ root of unity 
for $r$ the order of $g$. Write
$a(g)=\frac{1}{r}\sum a_i$ and  let $C(g)$ be the centralizer of $g\in G$.
Then one should have in the Grothendieck ring of varieties
$[Y]=\sum_{[g]} [X^g/C(g)][\A^{a(g)}]$.
Here $[g]$ runs through the conjugacy classes elements of $G$.
In particular we should have 
$A^i(Y)=\sum_{[g]} A^{i-a(g)}(X^g)^{C(g)}$.

Using  \cite{Go5},   Theorem \ref{mainthm} and the main result of 
\cite{dC-M2}  say that this is true
for the resolution
$\omega_n:S^{[n]}\to S^{(n)}$ of the symmetric power 
$S^{(n)}=S^n/\GG_n$.

\end{document}